\documentclass[a4paper,12pt,twoside]{article}
\usepackage{amsthm,amssymb,amsmath,amscd}
\usepackage{amsfonts,url}
\usepackage{latexsym}

\newtheorem{Thm}{Theorem}[section]
\newtheorem{lem}[Thm]{Lemma}

\newtheorem{de}[Thm]{Definition}

\newcommand{\q}{/\!/}
\newcommand{\N}{\mathbb{N}}

\newcommand{\T}{\mathbb{T}}
\newcommand{\Z}{\mathbb{Z}}
\newcommand{\C}{\mathbb{C}}
\newcommand{\G}{\mathcal{G}}
\newcommand{\X}{{\rm X}}

\newcommand{\fX}{\mathfrak{X}}
\def\SL{{\rm SL}}
\def\GL{{\rm GL}}

\def\fT{{\mathfrak T}}

\def\val{{\rm val}}

\def\SL{{\rm SL}}

\def\GL{{\rm GL}}
\def\PGL{{\rm PGL}}

\def\End{{\rm End}}

\def\Ind{{\rm Ind}}

\def\temp{{\rm temp}}

\def\Irr{{\rm Irr}}

\def\fs{{\mathfrak s}}

\def\fA{{\mathfrak A}}

\def\fK{{\mathfrak K}}

\def\fX{{\mathfrak X}}

\def\fY{{\mathfrak Y}}

\def\q{{/\!/}}

\title{$R$--\,groups and geometric structure in the representation theory of $\SL(N)$}
\author{Jamila Jawdat and Roger Plymen}
\date{}
\begin{document}

\maketitle \thispagestyle{empty}

\begin{abstract}Let $F$ be a nonarchimedean local field of characteristic zero and
let $G = \SL(N) = \SL(N,F)$.  This article is devoted to studying
the influence of the elliptic representations of $\SL(N)$ on the
$K$-theory.  We provide full arithmetic details. This study
reveals an intricate geometric structure. One point of interest is
that the $R$-group is realized as an isotropy group. Our results
illustrate, in a special case, part (3) of the recent conjecture
in \cite{ABP}.
\end{abstract}

\medskip


\bigskip

\section{Introduction}

Let $F$ be a nonarchimedean local field of characteristic zero and
let $G = \SL(N) = \SL(N,F)$.  This article is devoted to studying
subspaces of the tempered dual of $\SL(N)$ which have an
especially intricate geometric structure, and to computing, with
full arithmetic details, their $K$-theory. Our results illustrate,
in a special case, part (3) of the recent conjecture in
\cite{ABP}.

The subspaces of the tempered dual which are especially
interesting for us  contain \emph{elliptic} representations. A
tempered representation of $\SL(N)$ is \emph{elliptic} if its
Harish-Chandra character is not identically zero on the elliptic
set.

An element in the discrete series of $\SL(N)$ is an isolated point
in the tempered dual of $\SL(N)$ and contributes one generator to
$K_0$ of the reduced $C^*$-algebra of $\SL(N)$.

Now $\SL(N)$ admits elliptic representations which are not
discrete series: we investigate, with full arithmetic details, the
contribution of the elliptic representations of $\SL(N)$ to the
$K$-theory of the reduced $C^*$-algebra $\fA_N$ of $\SL(N)$.

According to \cite{Ply3}, $\fA_N$ is a $C^*$-direct sum of fixed
$C^*$-algebras.  Among these fixed algebras, we will focus on
those whose duals contain elliptic representations. Let $n$ be a
divisor of $N$ with $1 \leq n \leq N$ and suppose that the group
$\mathcal{U}_F$ of integer units admits a character of order $n$.
Then the relevant fixed algebras are of the form
\[
C(\T^n/\T,\mathfrak{K})^{\mathbb{Z}/n\mathbb{Z}} \subset \fA_N.\]
Here, $\mathfrak{K}$ is the $C^*$-algebra of compact operators on
standard Hilbert space, $\T^n/\T$ is the quotient of the compact
torus $\T^n$  via the diagonal action of $\T$. The compact group
$\T^n/\T$ arises as the maximal compact subgroup of the standard
maximal torus of the Langlands dual $\PGL(n,\C)$. We prove
(Theorem 3.1) that this fixed $C^*$-algebra is strongly Morita
equivalent to the crossed product
\[C(\T^n/\T) \rtimes \mathbb{Z}/n\mathbb{Z}.\]

The reduced $C^*$-algebra $\fA_N$ is liminal, and its primitive
ideal space is in canonical bijection with the tempered dual of
$\SL(N)$.  Transporting the Jacobson topology on the primitive
ideal space, we obtain a locally compact topology on the tempered
dual of $\SL(N)$, see \cite[3.1.1, 4.4.1, 18.3.2]{D}.

\smallskip
Let $\fT_n$ denote the $C^*$-dual of
$C(\T^n/\T,\fK)^{\mathbb{Z}/n\mathbb{Z}}$. Then $\fT_n$ is a
non-Hausdorff space, and has a very special structure as
topological space.  When $n$ is a prime number $\ell$, then
$\fT_{\ell}$ will contain multiple points. When $n$ is non-prime,
$\fT_n$ will contain not only multiple points, but also
\emph{multiple subspaces}. This crossed product $C^*$-algebra is a
noncommutative unital $C^*$-algebra which fits perfectly into the
framework of noncommutative geometry. In the tempered dual of
$\SL(N)$, there are connected compact non-Hausdorff spaces, laced
with  multiple subspaces, and  simply described by  crossed
 product $C^*$-algebras.

The $K$-theory of the fixed $C^*$-algebra is then given by the
$K$-theory of the crossed product $C^*$-algebra. To compute
(modulo torsion) the $K$-theory of this noncommutative
$C^*$-algebra, we apply the Chern character for discrete groups
\cite{BC}. This leads to the cohomology of the \emph{extended
quotient} $(\T^n/\T)\q (\mathbb{Z}/n\mathbb{Z})$. This in turn
leads to a problem in classical algebraic topology, namely the
determination of the cyclic invariants in the cohomology of the
$n$-torus.

The ordinary quotient will be denoted by $\fX(n)$:
 \[ \fX(n): =
(\T^n/\T)/(\Z/n\Z)\] This is a compact connected orbifold. Note
that $\fX(1) = pt$. The orbifold $\fX(n,k,\omega)$ which appears
in the following theorem is defined in section 4. The notation is
such that $\fX(n,n,1)$ is the ordinary quotient $\fX(n)$ and each
$\fX(n,1,\omega)$ is a point. The highest common factor of $n$ and
$k$ is denoted $(n,k)$.

\begin{Thm} The extended quotient $(\T^n/\T)\q (\Z/n\Z)$ is a disjoint union of compact connected orbifolds:
\[
(\T^n/\T)\q (\Z/n\Z)\, = \, \bigsqcup \; \fX(n,k,\omega)\] The
disjoint union is over all $1 \leq k \leq n$ and all $n/(k,n)$th
roots of unity $\omega$ in $\C$.
\end{Thm}

We apply the Chern character for discrete groups \cite{BC}, and
obtain
\begin{Thm} The K-theory groups $K_0$ and $K_1$ are given by
\[
K_0\, (C(\T^n/\T),\mathfrak{K})^{\Z/n\Z}\otimes_{\mathbb{Z}}\C
\,\simeq\, \bigoplus \; H^{ev}( \fX(n,k,\omega);\C)\]
\[K_1\, (C(\T^n/\T),\mathfrak{K})^{\Z/n\Z}\otimes_{\mathbb{Z}}\C \,\simeq\, \bigoplus \;
 H^{odd}(\fX(n,k,\omega);\C)\]
 The direct sums are over all $1 \leq k \leq n$ and all $n/(k,n)$th roots of unity $\omega$ in $\C$.
\end{Thm}

For the ordinary quotient  $\fX(n)$ we have the following explicit
formula (Theorems 6.1 and 6.3).   Let $H^{\bullet}: = H^{ev}
\oplus H^{odd}$ and let  $\phi$ denote the Euler totient.

\begin{Thm} Let $\fX(n)$ denote the ordinary quotient
$(\T^n/\T)/(\Z/n\Z)$. Then we have
\[dim_{\C}\,H^{\bullet}\, (\fX(n);\C) = \frac{1}{2n}\sum_{d|n,\, d\; odd}\phi(d)2^{n/d}.\]
\end{Thm}

  Theorem 1.1 lends itself to an interpretation in terms of representation theory.
 When $n = \ell$ a prime number, the elliptic representations of $\SL(\ell)$ are discussed in section 2.
 The extended quotient $(\T^{\ell}/\T)\q (\Z/ \ell Z)$ is the disjoint union of the ordinary quotient
 $\fX(\ell)$ and $\ell(\ell - 1)$ isolated points.
 We consider the canonical projection $\pi$ of the extended
 quotient onto the ordinary quotient:
  \[\pi: (\T^{\ell}/ \T)\q(\Z/\ell \Z) \longrightarrow \fX(\ell)\]
 The points $\tau_1, \ldots, \tau_{\ell}$ constructed in section
 2, are precisely the $\Z/\ell\Z$ fixed points in $\T^{\ell}/\T$.
 These are  $\ell$ points of reducibility, each of which admits $\ell$
 elliptic constituents.  Note also that, in the canonical
 projection $\pi$, the fibre $\pi^{-1}(\tau_j)$ of each point
 $\tau_j$ contains $\ell$ points. We may say that the extended
 quotient encodes, or provides a model of, reducibility. This is a very
 special case of the recent conjecture in \cite{ABP}.

 When $n$ is non-prime, we have points of reducibility, each of which admits elliptic constituents.
 In addition to the points of reducibility, there is a subspace of
reducibility. There are continua of $L$-packets. Theorem 1.2
describes the  contribution, modulo torsion, of all these
$L$-packets to $K_0$ and $K_1$.

 Let the infinitesimal character of the elliptic
representation $\epsilon$ be the cuspidal pair $(M,\sigma)$, where
$\sigma$ is an irreducible cuspidal representation of $M$ with
unitary central character. Then $\epsilon$ is a constituent of the
induced representation $i_{GM}(\sigma)$. Let $\fs$ be the point in
the Bernstein spectrum which contains the cuspidal pair
$(M,\sigma)$. To conform to the notation in \cite{ABP}, we will
write $E^{\fs}: = \T^n/\T,\, W^{\fs} = \Z/n\Z$. The standard
projection will be denoted
\[\pi^{\fs}:
 E^{\fs}\q W^{\fs} \to E^{\fs}/W^{\fs}.\]
 The space of tempered representations of $G$ determined by $\fs$ will be
 denoted  $\Irr^{\temp}(G)^{\fs}$, and the infinitesimal character will be denoted
 $inf.ch.$

 \begin{Thm}  There is a continuous bijection
  \[
 \mu^{\fs}: E^{\fs}\q W^{\fs} \longrightarrow
 \Irr^{\temp}(G)^{\fs}\] such that\[
 \pi^{\fs} = (inf.ch.) \circ \mu^{\fs}.\]
 This confirms, in a special case, part (3) of the conjecture in \cite{ABP}.
 \end{Thm}

 In section 2 of this article, we review elliptic representations
of the special linear algebraic group $\SL(N,F)$ over a $p$-adic
field $F$. Section 3 concerns fixed $C^*$-algebras and crossed
products. Section 4 computes the extended quotient $(\T^n/\T)\q
(\Z/n\Z)$. The formation of the $R$-groups is described in section
5. In section 6 we compute the cyclic invariants in the cohomology
of the $n$-torus.

 We would like to thank Paul Baum for several valuable discussions, Anne-Marie Aubert for her careful reading of the
 manuscript, Kuok Fai Chao and
 the referee for several constructive comments.

\section{The elliptic representations of $\SL(N)$}

Let $F$ be a nonarchimedean local field of characteristic zero.
Let $\textbf{G}$ be a connected reductive linear group over $F$.
Let $G=G(F)$ be the $F$-rational points of $\textbf{G}$. We say
that an element $x$ of $G$ is \emph{elliptic} if its centralizer
is compact modulo the center of $G$. We let $G^e$ denote the set
of regular elliptic elements of $G$.

Let $\mathcal{E}_{2}(G)$ denote  the set of equivalence classes of
irreducible discrete series representations of $G$, and denote by
$\mathcal{E}_{t}(G)$ be the set of equivalence classes of
irreducible tempered representations of $G$. Then
$\mathcal{E}_{2}(G) \subset \mathcal{E}_{t}(G)$.  If $\pi \in
\mathcal{E}_t(G)$, then we denote its character by $\Theta_{\pi}$.
Since $\Theta_{\pi}$ can be viewed as a locally integrable
function, we can consider its restriction to $G^e$, which we
denote by $\Theta^e_{\pi}$. We say that $\pi$ is elliptic if
$\Theta^e_{\pi} \neq 0$.  The set of elliptic representations
includes the discrete series.

Here is a classical example where elliptic representations occur
\cite{Assem}. We consider the group $\SL(\ell,F)$ with $\ell$ a
prime not equal to the residual characteristic of $F$. Let $K/F$
be a cyclic of order $\ell$ extension of $F$. The reciprocity law
in local class field theory is an isomorphism
\[
F^{\times}/N_{K/F}\,K^{\times} \cong \Gamma(K/F) = \mathbb{Z}/\ell
\mathbb{Z}\] where $\Gamma(K/F)$ is the Galois group of $K$ over
$F$. Let now $\mu_{\ell}(\C)$ be the group of $\ell$th roots of
unity in $\mathbb{C}$. A choice of isomorphism $\Z/\ell\Z \cong
\mu_{\ell}(\C)$ then produces a character $\kappa$ of $F^{\times}$
of order $\ell$ as follows:
\[
\kappa : F^{\times} \to F^{\times}/N_{K/F}\,K^{\times} \cong
\mathbb{Z}/\ell\mathbb{Z} \cong \mu_{\ell}(\mathbb{C})\]

Let $B$ be the standard Borel subgroup of $\SL(\ell)$, let $T$ be
the standard maximal torus, and let $B = T\cdot N$ be its Levi
decomposition. Let $\tau$ be the character of $T$ defined by
\[
\tau: = 1 \otimes \kappa \otimes \cdots \otimes \kappa^{\ell -
1}\] and let\[\pi(\tau):=\, \Ind^G_B(\tau \otimes 1)\]be the
unitarily induced representation of $\SL(\ell)$.

Now $\pi(\tau)$ is a representation in the minimal unitary
principal series of $\SL(\ell)$. It has $\ell$ distinct
irreducible elliptic components and the Galois group $\Gamma(K/F)$
acts simply transitively on the set of irreducible components. The
set of irreducible components of $\pi(\tau)$ is an $L$-packet.

Let\[ \pi(\tau) = \pi_1 \oplus \cdots \oplus \pi_{\ell}\] be the
$\ell$ components of $\pi(\tau)$. The character $\Theta$ of
$\pi(\tau)$, as character of a principal series representation,
\emph{vanishes on the elliptic set}.  The character $\Theta_1$ of
$\pi_1$ on the elliptic set is therefore \emph{cancelled out} by
the sum $\Theta_2 + \cdots + \Theta_{\ell}$ of the characters of
the relatives $\pi_2, \ldots, \pi_{\ell}$ of $\pi_1$.

Let $\omega$ denote an $\ell$th root of unity in $\C$. All the
$\ell$th roots are allowed, including $\omega = 1$.  In the
definition of $\tau$, we now replace $\kappa$ by $\kappa \otimes
\omega^{\val}$. This will create $\ell$ characters, which we will
denote by $\tau_1, \ldots, \tau_{\ell}$, where $\tau_1 = \tau$.
For each of these characters, the $R$-group is given as follows:
\[
R(\tau_j) = \Z/\ell \Z
\] for all $1 \leq j \leq \ell$, and the
induced representation $\pi(\tau_j)$ admits $\ell$ elliptic
constituents.


If $P=MU$ is a standard parabolic subgroup of $G$ then
$i_{GM}(\sigma)$ will denote the induced representation
$\Ind_{MU}^G(\sigma \otimes 1)$ (normalized induction).  The
$R$-group attached to $\sigma$ will be denoted $R(\sigma)$.

Let $P=MU$ be the standard parabolic subgroup of $G: = \SL(N,F)$
described as follows.  Let $N = mn$, let $\tilde{M}$ be the Levi
subgroup $\GL(m)^n \subset \GL(N,F)$  and let $M=\tilde{M} \cap
\SL(N,F)$.

We will use the framework, notation and main result in \cite{Gol}.
Let $\sigma \in \mathcal{E}_2(M)$ and let $\pi_{\sigma} \in
\mathcal{E}_2(\widetilde{M})$ with $\pi_{\sigma}|M \supset
\sigma$.  Let $W(M): = N_G(M)/M$ denote the Weyl group of $M$, so
that $W(M)$ is the symmetric group on $n$ letters.
Let\[\overline{L}(\pi_{\sigma}): = \{\eta \in \widehat{F^{\times}}
| \pi_{\sigma} \otimes \eta \simeq w\pi_{\sigma} \quad \textrm{for
some} \quad w \in W\}\]

\[X(\pi_{\sigma}): =
\{\eta \in \widehat{F^{\times}} | \pi_{\sigma} \otimes \eta \simeq
\pi_{\sigma}\}\]

By \cite[Theorem 2.4]{Gol}, the $R$-group of $\sigma$ is given by
\[
R(\sigma)  \simeq \overline{L}(\pi_{\sigma})/X(\pi_{\sigma}).\]

We follow \cite[Theorem 3.4]{Gol}. Let $\eta$ be a smooth
character of $F^{\times}$ such that $\eta^n \in X(\pi_1)$ and
$\eta^j \notin X(\pi_1)$ for $1 \leq j \leq n-1$. Set
\begin{equation}\label{eq:1}
\pi_{\sigma} \simeq \pi_1 \otimes \eta\pi_1 \otimes \eta^2\pi_1
\otimes \cdots \otimes \eta^{n-1}\pi_1, \quad \quad \pi_{\sigma}|M
\supset \sigma\end{equation}with $\pi_1 \in
\mathcal{E}_2(\GL(m))$, $\eta\pi_1: = (\eta \circ \det) \otimes
\pi_1$. Then we have
\[\overline{L}(\pi_{\sigma})/X(\pi_{\sigma}) = <\eta>\] and so
$R(\sigma) \simeq \Z/n\Z$.  The elliptic representations
are the constituents of $i_{GM}(\sigma)$ with $\pi_{\sigma}$ as in
equation (\ref{eq:1}).

\section{Fixed algebras and crossed products}
 Let $M$ denote the Levi subgroup which occurs in section 2. Denote by $\Psi^1(M)$
the group of unramified unitary characters of $M$. Now $M \subset
\SL(N,F)$ comprises blocks $x_1,\cdots, x_n$ with $x_i \in
\GL(m,F)$ and $\prod \det(x_i)=1$. Each unramified unitary
character $\psi \in \Psi^1(M)$ can be expressed as follows,\[\psi:
diag(x_1, \cdots, x_n) \to \prod_{j = 1}^n z_j^{\val(\det x_j)}\]
with $z_1,z_2, \cdots ,z_n \in \mathbb{T}$, i.e. $|z_i|=1$. Such
unramified unitary characters $\psi$ correspond to coordinates
$(z_1:z_2: \cdots :z_n)$ with each $z_i \in \mathbb{T}$. Since
 \[\prod_{i = 1}^n(zz_i)^{\val(\det x_i)} = \prod_{i=1}^n z_i^{\val(\det x_i)}\]
  we have \emph{homogeneous} coordinates.
We  have the isomorphism
\[
\Psi^1(M) \cong \{(z_1:z_2: \cdots :z_n): \, |z_i|=1,1 \leq i \leq
n\} = \T^n/\T.
\]

If $M$ is the standard maximal torus $T$ of $\SL(N)$ then
$\Psi^1(T)$ is the maximal \emph{compact} torus in the dual torus
\[ T^{\vee} \subset G^{\vee} = \PGL(N,\C)\] where $G^{\vee}$ is
the Langlands dual group.

Let $\sigma, \, \pi_{\sigma}, \, \pi_1$  be as in equation
(\ref{eq:1}). Let $g$ be the order of the group of unramified
characters $\chi$ of $F^{\times}$ such that $(\chi \circ \det)
\otimes \pi_1 \simeq \pi_1$. Now let \[E: = \{\psi\otimes \sigma:
\psi \in \Psi^1(M)\}.\] The base point $\sigma \in E$ determines a
homeomorpism
\[E \simeq \T^n/\T, \quad (z_1^{\val \circ \det} \otimes \cdots \otimes
z_n^{\val \circ \det})\otimes \sigma \mapsto (z_1^g: \cdots :
z_n^g).\]

 From this point onwards, we will require that the \emph{restriction of
$\eta$ to the group $\mathcal{U}_F$ of integer units is of order
$n$.} Let $W(M)$ denote the Weyl group of $M$ and let $W(M,E)$ be
the subgroup of $W(M)$ which leaves $E$ globally invariant. Then
we have $W(M,E) = W(\sigma) = R(\sigma) = \Z/n\Z$.

Let $\mathfrak{K} = \mathfrak{K}(H)$ denote the $C^*$-algebra of
compact operators on the standard Hilbert space $H$. Let
$\mathfrak{a}(w,\lambda)$ denote normalized intertwining
operators. The fixed $C^*$-algebra $C(E,\mathfrak{K})^{W(M,E)}$ is
given by
\[\{f \in C(E,\mathfrak{K}) |
f(w\lambda)=\mathfrak{a}(w,\lambda \tau)f(\lambda)\mathfrak{a}(w,
\lambda \tau)^{-1}, w \in W(M,E)\}.\] This fixed $C^*$-algebra is
a $C^*$-direct summand of the reduced $C^*$-algebra $\fA_N$ of
$\SL(N)$, see \cite{Ply3}.

\begin{Thm} Let  $G= \SL(N,F)$, and $M$ be a Levi subgroup
consisting of $n$ blocks of the same size $m$. Let $\sigma \in
\mathcal{E}_2(M)$. Assume that the induced representation
$i_{GM}(\sigma)$ has  elliptic constituents, then the fixed
$C^*$-algebra $C(E,\mathfrak{K})^{W(M,E)}$ is strongly Morita
equivalent to the crossed product $C^*$-algebra $C(E) \rtimes
\mathbb{Z}/n\Z$.
\end{Thm}

\begin{proof}  For the commuting algebra of $i_{MG}(\sigma)$, we
have \cite{Sil}:
\[
\End_G((i_{MG}(\sigma)) = \mathbb{C}[R(\sigma)].\]

Let $w_0$ be a generator of $R(\sigma)$, then the normalized
intertwining operator $\mathfrak{a}(w_0,\sigma)$ is a unitary
operator of order $n$.
By the spectral theorem for unitary operators, we have
\[\mathfrak{a}(w_0,\sigma)= \sum_{j=0}^{n-1} ~\omega^{j}\, E_{j} \]where  $\omega = \exp(2\pi i/n)$ and
$E_j$ are the projections onto the irreducible subspaces of the
induced representation $i_{MG}(\sigma)$. The unitary
representation\[R(\sigma) \to U(H), \quad w \mapsto
\mathfrak{a}(w,\sigma)\] contains each character of $R(\sigma)$
countably many times. Therefore condition (***) in \cite[p.
301]{PL} is satisfied. The condition (**) in \cite[p. 300]{PL} is
trivially satisfied since $W(\sigma) = R(\sigma)$.

We have $W(\sigma) = \Z/n\Z$. Then a subgroup $W(\rho)$ of order
$d$ is given by $W(\rho) = k\Z \mod n$ with $dk = n$. In that
case, we have
\[
\mathfrak{a}(w_0,\sigma)|_{W(\rho)}= \sum_{j = 0}^{n-1}
    \omega^{kj}E_{j}.\]
We compare the two unitary representations:
\[\phi_1: W(\rho)\to U(H), \quad  w \mapsto
\mathfrak{a}(w,\sigma)|_{W(\rho)}\]
\[\phi_2: W(\rho) \to U(H),
\quad w \mapsto \mathfrak{a}(w,\rho).\] Each representation
contains every character of $W(\rho)$. They are
\emph{quasi-equivalent} as in \cite{PL}. Choose an increasing
sequence $(e_n)$ of finite-rank projections in $\mathcal{L}(H)$
which converge strongly to $I$ and commute with each projection
$E_j$. The compressions of $\phi_1, \phi_2$ to $e_{n}H$ remain
quasi-equivalent. Condition (*) in \cite[p. 299]{PL} is satisfied.

All three conditions of \cite[Theorem 2.13]{PL} are satisfied. We
therefore have a strong Morita equivalence
\[(C(E)\otimes \mathfrak{K})^{W(M,E)} \simeq C(E) \rtimes R(\sigma)=
\C(E) \rtimes \Z/n\Z. \]\end{proof}

We will need a special case of the Chern character for discrete
groups \cite{BC}.
 \begin{Thm} We have an isomorphism
 \[K_i(C(E) \rtimes \Z/n\Z) \otimes_{\mathbb{Z}} \mathbb{C} \cong
 \bigoplus_{j\in \N} H^{2j+i}(E\q(\Z/n\Z) ;\C)\]with $i
 = 0,1$, where $E\q(\Z/n\Z)$ denotes the extended quotient of
 $E$ by $\Z/n\Z$.
\end{Thm}

When $N$ is a prime number $\ell$, this result already appeared in
\cite{Ply2, PL}.

\section{The formation of the fixed sets}

Extended quotients were introduced by Baum and Connes \cite{BC} in
the context of the Chern character for discrete groups. Extended
quotients were used in \cite{Ply1, Ply2} in the context of the
reduced group $C^*$-algebras of $\GL(N)$ and $\SL(\ell)$ where
$\ell$ is prime.   The results in this section extend results in
\cite{Ply2, PL}.

\begin{de}  Let $\rm X$ be a  compact
Hausdorff topological space. Let $\Gamma$ be a finite
\emph{abelian} group acting on $\rm X$ by a (left) continuous
action. Let
\[
\widetilde{\rm{X}}=\{(x,\gamma) \in {\rm{X}} \times \Gamma :
\gamma x=x\}\] with the group action on  $\widetilde {\rm{X}}$
given by
\[g \cdot (x,\gamma)=(gx,\gamma)\] for $g \in
\Gamma$.  Then the
 \emph{extended quotient} is given by \[ \rm X\q\Gamma: = \widetilde{{\rm X}}/\Gamma = \bigsqcup_{\gamma \in \;\Gamma}
 \rm X^{\gamma}/ \Gamma\]where $\rm X^{\gamma}$ is the $\gamma$-fixed set. \end{de}

The extended quotient will always contain the ordinary quotient.
The standard projection $\pi: X\q \Gamma \to X/\Gamma$ is induced
by the map $(x,\gamma) \mapsto x$.  We note the following
elementary fact, which will be useful later (in Lemma 5.2): let $y
= \Gamma x$ be a point in $X/\Gamma$. Then the cardinality of the
pre-image $\pi^{-1}y$ is equal to the order of the isotropy group
$\Gamma_x$: \[ |\pi^{-1}y| = |\Gamma_x|.\]

We will write $\X = E = \T^n/\T$, where $\T$ acts diagonally on
$\T^n$, i.e.
\[t(t_1,t_2, \cdots ,t_n)= (tt_1,tt_2, \ldots ,tt_n), \quad \quad t, t_i \in \T.\] We
have the action of the finite group $\Gamma = \Z/n\Z$ on $\T^n/\T$
given by cyclic permutation. The two actions of $\T$ and of
$\Z/n\Z$ on $\T^n$ commute.
We will write $(k,n)$ for the highest common factor of $k$ and
$n$.

\begin{Thm} The extended quotient $(\T^n/\T)\q (\Z/n\Z)$ is a disjoint union of
compact connected orbifolds:
\[(\T^n/\T)\q (\Z/n\Z) \;\simeq \; \bigsqcup_{1\leq k \leq n, \;\omega^{n/(k,n)}=1}
\fX(n,k,\omega)\] Here, $\omega$ is a $n/(k,n)$th root of unity in
$\C$.
\end{Thm}

\begin{proof} Let $\gamma$ be the standard $n$-cycle defined by $\gamma(i) = i + 1 \mod
n$. Then $\gamma^k$ is the product of $n/d$ cycles of order $d =
n/(n,k)$. Let $\omega$ be a $d$th root of unity in $\C$.  All
$d$th roots of unity are allowed, including $\omega = 1$.  The
element $t(\omega) = t(\omega; z_1, \ldots,z_n) \in \T^n$ is
defined by imposing the following relations:
\[
z_{i + k} = \omega^{-1} z_i\]all suffices $\mod n$. This condition
allows $n/d$ of the complex numbers $z_1, \ldots, z_n$ to vary
freely, subject only to the condition that each $z_j$ has modulus
$1$. The crucial point is that
\[\gamma^k \cdot t(\omega) = \omega t(\omega)\] Then $\omega$
determines a $\gamma^k$-fixed set in $\T^n/\T$, namely the set
$\fY(n,k,\omega)$ of all cosets $t(\omega)\cdot \T$.
 The set  $\fY(n,k,\omega)$ is an
$(n/d-1)$-dimensional subspace of fixed points.

Note that $\fY(n,k,\omega)$, as a coset of the closed subgroup
$\fY(n,k,1)$ in the compact Lie group $E$, is homeomorphic (by
translation in $E$) to $\fY(n,k,1)$.  The translation is by the
element $t(\omega: 1, \ldots,1)$. If $\omega_1, \omega_2$ are
distinct $d$th roots of unity, then $\fY(n,k,\omega_1),
\fY(n,k,\omega_2)$ are disjoint.

We define the quotient space \[\fX(n,k,\omega): =
\fY(n,k,\omega)/(\Z/n\Z)\]and apply definition 4.1.
\end{proof}

When $k = n$, we must have $\omega = 1$. In that case, the
orbifold is the ordinary quotient: $\fX(n,n,1) = \fX(n)$.

Let $(n,k) = 1$.  The number of such $k$ in $1 \leq k \leq n$ is
$\phi(n)$. In this case,  $\omega$ is an $n$th root of unity and
$\fX(n,k,\omega)$ is a point.  There are $n$ such roots of unity
in $\C$.  Therefore, the extended quotient $(\T^n/\T)\q (\Z/n\Z)$
always contains $\phi(n)n$ isolated points.

Theorem 1.1 is a consequence of Theorems 3.1, 3.2 and 4.2. If, in
Theorem 1.1, we take $n$ to be a prime number $\ell$, then we
recover the following result in \cite[p. 30]{Ply2}: the extended
quotient $(\T^{\ell}/\T)\q(\Z/\ell\Z)$ is the disjoint union of
the ordinary quotient $\fX(\ell)$ and $(\ell - 1)\ell$ points.

\section{The formation of the $R$-groups}

We continue with the notation of section 3. Let
$\sigma,\,\pi_{\sigma}, \, \pi_1, \eta$ be as in equation
(\ref{eq:1}). The $n$-tuple $t: = (z_1, \ldots,z_n)\in \T^n$
determines an element $[t] \in E$. We can interpret $[t]$ as the
unramified character
\[
\chi_t: = (z_1^{\val\circ\det}, \ldots, z_n^{\val\circ\det})\]

Let $\Gamma = \Z/n\Z$, and let $\Gamma_{[t]}$ denote the isotropy
subgroup of $\Gamma$.
\begin{lem} The isotropy subgroup $\Gamma_{[t]}$ is isomorphic to
the $R$-group of $\chi_t\otimes \sigma$:
\[
\Gamma_{[t]} \simeq R(\chi_t \otimes \sigma)\]
\end{lem}
\begin{proof} Let the order of $\Gamma_{[t]}$ be $d$. Then $d$ is a
divisor of $n$.  Let $\gamma$ be a generator of $\Gamma_{[t]}$.
Then $\gamma$ is a product of $n/d$ disjoint $d$-cycles, as in
section 4.  We must have $t = t(\omega)$ with $\omega$ a $d$th
root of unity in $\C$. Note that $\gamma \cdot t(\omega) = \omega
t(\omega)$. Then we have
\begin{eqnarray*} R(\chi_t\otimes \sigma) &=&
\overline{L}(\chi_t\otimes\pi_{\sigma})/X(\chi_t\otimes\pi_{\sigma})
\nonumber\\
 &=& \{\alpha \in \widehat{F^{\times}}: w \pi_{\sigma}
\simeq \pi_{\sigma} \otimes
\alpha\;  \textrm{for some $w$ in}\; W\}/X(\chi_t\otimes\pi_{\sigma})    \nonumber\\
&=& <\omega^{\val\circ\det}\otimes \eta^{n/d}>     \nonumber\\
&=& \Z/d \Z \nonumber\\
&=& \Gamma_{[t]}\nonumber\\
\end{eqnarray*}
since, modulo $X(\chi_t\otimes\pi_{\sigma})$, the character
$\eta^{n/d}$ has order $d$.
\end{proof}
\begin{lem}In the standard projection $p: E\q \Gamma \to E/\Gamma$,
the cardinality of the fibre of $[t]$ is the order of the
$R$-group of $\chi_t\otimes \sigma$.\end{lem}
\begin{proof} This follows from Lemma 5.1.
\end{proof}

We will assume that $\sigma$ is a \emph{cuspidal} representation
of $M$ with unitary central character. Let $\fs$ be the point in
the Bernstein spectrum of $\SL(N)$ which contains the cuspidal
pair $(M,\sigma)$. To conform to the notation in \cite{ABP}, we
will write $E^{\fs}: = \T^n/\T,\, W^{\fs} = \Z/n\Z$. The standard
projection will be denoted
\[\pi^{\fs}:
 E^{\fs}\q W^{\fs} \to E^{\fs}/W^{\fs}.\]
 The space of tempered representations of $G$ determined by $\fs$ will be
 denoted  $\Irr^{\temp}(G)^{\fs}$, and the infinitesimal character will be denoted
 $inf.ch.$

\begin{Thm} We have a commutative diagram:

\[
\begin{CD}
E\q W^{\fs} @>\mu^{\fs}>> \Irr^{\temp}(G)^{\fs}\\
@V\pi^{\fs} VV @VVinf.ch. V\\
E/W^{\fs} @>>> E/W^{\fs}
\end{CD}
\]
in which the map $\mu^{\fs}$ is a continuous bijection. This
confirms, in a special case, part (3) of the conjecture in
\cite{ABP}.
\end{Thm}

 \begin{proof} We have
 \[
 \C[R(\sigma)] \simeq \End_G(i_{GM}(\sigma))
 \]
 This implies that the characters of the cyclic group $R(\sigma)$
 parametrize the irreducible constituents of $i_{GM}(\sigma)$.  This
 leads to a labelling of the irreducible constituents of
 $i_{GM}(\sigma)$, which we will write as $i_{GM}(\sigma:r)$ with
 $0 \leq r < n$.

 The map $\mu^{\fs}$ is defined as follows:
 \[
 \mu^{\fs}: (t,\gamma^{rd}) \mapsto i_{GM}(\chi_t\otimes\sigma:r)\]
We now apply Lemma 5.2.

  Theorem 3.2 in \cite{Ply3} relates the natural topology on the
Harish-Chandra parameter space to the Jacobson topology on the
tempered dual of a reductive $p$-adic group. As a consequence, the
map $\mu^{\fs}$ is continuous.
\end{proof}

\section{Cyclic invariants}

We will consider the map \[ \alpha : \T^n \to (\T^n /\T) \times
\T, \quad \quad (t_1, \ldots, t_n) \to ((t_1: \ldots :t_n), t_1
t_2 \cdots t_n)\]where $(t_1: \ldots : t_n)$ is the image of
$(t_1, \ldots,t_n)$ via the map $\T^n \to \T^n/\T$. The map
$\alpha$ is a homomorphism of Lie groups. The kernel of this map
is
\[
\G_n: = \{\omega I_n: \omega^n = 1\}.\] We therefore have the
isomorphism of compact connected Lie groups:
\begin{eqnarray}
\T^n/\G_n \; \cong \; (\T^n /\T) \times \mathbb{T}
\end{eqnarray}
This isomorphism is equivariant with respect to the $\Z/n\Z$-
action, and we infer that
\begin{eqnarray}
(\T^n/\G_n)/(\Z/n\Z) \; \cong \; (\T^n /\T)/(\Z/n\Z) \times
\mathbb{T}
\end{eqnarray}

\begin{Thm} Let
 $H^{\bullet}(-;\C)$ denote the total cohomology group. We have
\[
\dim_{\C}\,H^{\bullet}(\fX(n);\C) =\, \frac{1}{2} \cdot
\dim_{\C}\,H^{\bullet}(\mathbb{T}^n;\C)^{\Z/n\Z}.\]
\end{Thm}
\begin{proof} The cohomology of the orbit space is
given by the fixed set of the cohomology of the original space
\cite[Corollary 2.3, p.38]{Bor}.  We have
\begin{eqnarray}H^j(\T^n/\G_n;\C)& \cong &
H^j(\T^n;\C)^{\G_n}\nonumber\\
& \cong & H^j(\T^n;\C)
\end{eqnarray}
since the action of $\G_n$ on $\T^n$ is homotopic to the identity.
We spell this out. Let $z: = (z_1, \ldots,z_n)$ and define $H(z,t)
= \omega^t \cdot z = (\omega^t z_1, \ldots,\omega^t z_n)$. Then
$H(z,0) = z, \; H(z,1) = \omega \cdot z$. Also, $H$ is equivariant
with respect to the permutation action of $\Z/n\Z$. That is to
say, if $\epsilon \in \Z/n\Z$ then $H(\epsilon\cdot z,t) =
\epsilon\cdot H(z,t)$. This allows us to proceed as follows:
\begin{eqnarray}H^j(\T^n;\C)^{\Z/n\Z}& \cong &
H^j(\T^n/\G_n;\C)^{\Z/n\Z}\nonumber\\
& \cong & H^j((\T^n/\T) \times \T;\C)^{\Z/n\Z}\nonumber\\
& \cong & H^j((\T^n/\T)/(\Z/n\Z)\times \T;\C)\label{coh}
\end{eqnarray}
We apply the Kunneth theorem in cohomology (there is no torsion):
\[(H^j(\T^n;\C))^{\Z/n\Z} \cong  H^j(\fX(n);\C) \oplus
H^{j-1}(\fX(n);\C)\quad \quad \textrm{with}\; 0 < j \leq n\]
\[(H^n(\T^n;\C))^{\Z/n\Z} \simeq H^{n-1}(\fX(n);\C), \quad H^0(\T^n;\C)^{\Z/n\Z} \cong  H^0(\fX(n);\C) \simeq \C\]

\[ H^{ev}(\T^n;\C)^{\Z/n\Z} = H^{\bullet}(\fX(n);\C),\quad H^{odd}(\T^n;\C)^{\Z/n\Z} = H^{\bullet}(\fX(n);\C)\]
\end{proof}

We now have to find the cyclic invariants in
$H^{\bullet}(\T^n;\C)$. The cohomology ring $H^{\bullet}(\T^n,
\C)$ is the exterior algebra $\bigwedge V$ of a complex
$n$-dimensional vector space $V$, as can be seen by considering
differential forms $d\theta_1 \wedge \cdots \wedge d\theta_r$. The
vector space $V$ admits a basis $\alpha_1 = d\theta_1,
\ldots,\alpha_n = d\theta_n$. The action of $\Z/n\Z$ on $\bigwedge
V$ is induced by permuting the elements $\alpha_1,
\ldots,\alpha_n$, i.e. by the regular representation $\rho$ of the
cyclic group $\Z/n\Z$. This representation of $\Z/n\Z$ on
$\bigwedge V$ will be denoted $\bigwedge \rho$. The dimension of
the space of cyclic invariants in $H^{\bullet}(\T^n,\C)$ is equal
to the multiplicity of the unit representation $1$ in $\bigwedge
\rho$. To determine this, we use the theory of group characters.

\begin{lem} The dimension of the the subspace of cyclic invariants is given by
\[(\chi_{\bigwedge \rho}, 1) =
\frac{1}{n}(\chi_{\bigwedge \rho}(0) + \chi_{\bigwedge \rho}(1) +
\cdots + \chi_{\bigwedge \rho}(n-1)).\]
\end{lem}
\begin{proof} This is a standard result in the theory of group
characters \cite{S}.
\end{proof}

\begin{Thm} The dimension of the space of cyclic invariants in $H^{\bullet}(\T^n,\C)$
is given by the formula \[g(n): = \frac{1}{n}\sum_{d|n, d\,
\textrm{odd}}\phi(d)2^{n/d}\]
\end{Thm}

\begin{proof} We note first that

\[\chi_{\bigwedge \rho}(0) = \textrm{Trace}\;1_{\bigwedge V} = \textrm{dim}_{\C}\;\bigwedge V = 2^n.\]

To evaluate the remaining terms, we need to recall the definition
of the elementary symmetric functions $e_j$:
\[
\prod_{j=1}^{n}(\lambda - \alpha_j) = \lambda^n -\lambda^{n-1}e_1
+ \lambda^{n-2}e_2 - \cdots + (-1)^ne_n.\] When we need to mark
the dependence on $\alpha_1,\ldots,\alpha_n$ we will write $e_j =
e_j(\alpha_1,\ldots,\alpha_n)$. Set $\alpha_j = \omega^{j-1},
\omega = \exp(2\pi i/n)$. Then we get
\[\lambda^n - 1 = \prod_{j=1}^{n}(\lambda - \alpha_j) = \lambda^n -\lambda^{n-1}e_1
+ \lambda^{n-2}e_2 - \cdots + (-1)^ne_n.\]

Let $d|n$, let $\zeta$ be a \emph{primitive} $d$th root of unity.
Let $\alpha_j = \zeta^{j-1}$. We have
\begin{eqnarray}(\lambda^d-1)^{n/d} =(\lambda^d-1) \cdots (\lambda^d-1) =
\prod_{j=1}^n (\lambda - \alpha_j) \end{eqnarray} Set $\lambda =
-1$. If $d$ is even, we obtain
\begin{eqnarray}0= 1 + e_1(1,\zeta,\zeta^2,\ldots) +
e_2(1,\zeta,\zeta^2,\ldots) + \cdots +
e_{n}(1,\zeta,\zeta^2,\ldots)\end{eqnarray} If $d$ is odd, we
obtain
\begin{eqnarray}
2^{n/d} =   1 + e_1(1,\zeta,\zeta^2,\ldots) +
e_2(1,\zeta,\zeta^2,\ldots) + \cdots +
e_n(1,\zeta,\zeta^2,\ldots)\end{eqnarray}

We observe that the regular representation $\rho$ of the cyclic
group $\Z/n\Z$ is a direct sum of the characters $m \mapsto
\omega^{rm}$ with $0 \leq r \leq n $. This direct sum
decomposition allows us to choose a basis $v_1, \ldots,v_n$ in $V$
such that the representation $\bigwedge \rho$ is diagonalized by
the wedge products $v_{j_1}\wedge \cdots \wedge v_{j_l}$. This in
turn allows us to compute the character of $\bigwedge \rho$ in
terms of the elementary symmetric functions $e_1, \ldots,e_n$.

With $\zeta = \omega^r$ as above, we have
\[
\chi_{\bigwedge \rho}(r) = 1 + e_1(1,\zeta, \zeta^2, \ldots) +
e_2(1,\zeta,\zeta^2,\ldots) + \cdots + e_n(1,\zeta,\zeta^2,
\ldots)\]

We now sum the values of the character $\chi_{\bigwedge \rho}$.
Let $d: = n/(r,n)$. Then $\zeta$ is a primitive $d$th root of
unity. If $d$ is even then $\chi_{\bigwedge \rho}(r) = 0$. If $d$
is odd, then $\chi_{\bigwedge \rho}(r) = 2^{n/d}$. There are
$\phi(d)$ such terms. So we have
\begin{eqnarray} \chi_{\bigwedge \rho}(0)+ \chi_{\bigwedge \rho}(1) + \cdots +\chi_{\bigwedge \rho}(n-1)
= \sum_{d|n,\, d\, \textrm{odd}}\phi(d)2^{n/d}
\end{eqnarray}
We now apply Lemma 6.2.\end{proof}

The sequence $n \mapsto g(n)/2$, $n = 1,2,3,4, \ldots$, is
\[1,1,2,2,4,6,10,16,30,52,94,172,316,586,1096,2048,3856, 7286,\ldots\]
 as in
\url{www.research.att.com/~njas/sequences/A000016}. Thanks to
Kasper Andersen for alerting us to this web site.

J. Jawdat, Department of Mathematics, Zarqa Private
University, Jordan.\\

Email: jjawdat@zpu.edu.jo\\

R. Plymen, School of Mathematics, Alan Turing Building,
Manchester University, Manchester M13 9PL, England.\\

Email: plymen@manchester.ac.uk\\

\end{document}